\documentclass{amsart}
\usepackage{amssymb}
\usepackage[mathscr]{eucal}

\newtheorem{theorem}{Theorem}[section]
\newtheorem{corollary}[theorem]{Corollary}
\newtheorem{lemma}[theorem]{Lemma}
\newtheorem{proposition}[theorem]{Proposition}

\theoremstyle{definition}
\newtheorem{definition}[theorem]{Definition}

\newtheorem{example}[theorem]{Example}

\theoremstyle{remark}
\newtheorem{remark}[theorem]{Remark}

\begin{document}

\title{An operator Arzel\`a-Ascoli theorem}

\author{wei wu}
\address{Department of Mathematics, East China Normal University, Shanghai 200062, P.R. China}
\email{wwu@math.ecnu.edu.cn}
\curraddr{Department of Mathematics, University of California, Berkeley, CA 94720-3840}
\email{wwu@math.berkeley.edu}

\begin{abstract}
We generalize the Arzel\`a-Ascoli theorem to the setting of matrix order unit spaces, extending the work of Antonescu-Christensen on 
unital $C^\ast$-algebras. This gives an affirmative answer to a question of Antonescu and Christensen.
\end{abstract}

\maketitle

\section{Introduction}\label{se:1}

Motivated by the observation that for a compact spin Riemannian manifold one can recover its smooth structure, its Riemannian metric, 
and much else, directly from its standard Dirac operator, Connes pointed out that from a spectral triple $(\mathcal{A}, \mathcal{H}, D)$ one 
obtains a metric on the state space $\mathcal{S}(\mathcal{A})$ of the unital $C^\ast$-algebra $\mathcal{A}$ by the formula
\[\rho_D(\varphi,\psi)=\sup\{|\varphi(a)-\psi(a)|: \|[D, a]\|\le 1, a\in \mathcal{A}\},\]
if $\{a\in \mathcal{A}; \|[D, a]\|\le 1\}/(\mathbb{C}1)$ is bounded\cite{co1}.

A natural question is that when the topology on $\mathcal{S}(\mathcal{A})$ determined by $\rho_D$ coincides with the $w^\ast$-topology. Rieffel has 
studied it in a more general situation in which $\mathcal{A}$ is just an order unit space and $\|[D, \cdot]\|$ is replaced by a Lipschitz 
seminorm on $\mathcal{A}$\cite{ri1, ri2, ri3}. This and certain statements in the high energy physics and string theory, concerning 
non-commutative spaces that converged to other spaces, led him to the concept of compact quantum metric spaces\cite{ri4, ri5}.

For a discrete group $G$, which is of rapid decay with respect to some length function, Antonescu and Christensen got a metric on the 
state space $\mathcal{S}(C_r^\ast(G))$ of the reduced group $C^\ast$-algebra $C_r^\ast(G)$ which is finite for all pairs, is bounded on 
$\mathcal{S}(C_r^\ast(G))$, and generates the $w^\ast$-topology on $\mathcal{S}(C_r^\ast(G))$\cite{anch}. This way of obtaining a metric 
from higher derivatives inspired them to discuss metric spaces without a smooth structure. They believe that any norm compact balanced 
convex subset of a unital $C^\ast$-algebra $\mathcal{A}$ which separates the states on $\mathcal{A}$ contains much information needed, and call the subset 
a {\it metric set} of $\mathcal{A}$. In particular, they showed that it works well with respect to a translation of the classical Arzel\`a-Ascoli theorem 
into a non-commutative language\cite{anch}.

Most of interesting constructions in view of Lipschitz seminorms on $C^\ast$-algebras, such as those from Dirac operators, or those in 
\cite{ri1}, also provide in a natural way seminorms on 
all the matrix algebras over the algebras. Rieffel suggested that some ``matrix Lipschitz seminorm" in analogy with the matrix 
norms of \cite{ef} will be of importance\cite{ri2}. In \cite{wu1,wu2}, we developed a version of it on the matrix order unit 
spaces. It has many nice properties\cite{wu1, wu2, wu3}. In \cite{anch}, Antonescu and Christensen asked if their result on the 
non-commutative version of the Arzel\`a-Ascoli theorem is valid in a wider generality like operator systems. The main goal 
of this paper is to give an affirmative answer to their question both at the ``matrix" level and at the ``function" level.

This paper is organized as follows. We begin in Section \ref{se:2} with a discussion of the notions. Because we need corresponding 
equicontinuity of mappings on the framework of operator spaces, we discuss the continuous matrix mappings on matrix metric spaces in 
Section \ref{se:3}. We introduce the concept of matrix metric sets of a matrix order unit space $(A, 1)$ which is closely related 
to matrix metrics on the matrix state space $\mathcal{CS}(A)$ generating the BW-topology. This is done in Section \ref{se:4}. We prove 
our main results (Theorem \ref{th:54} and Corollary \ref{co:55}) on relative compactness in Section \ref{se:5}.

\section{Preliminaries}\label{se:2}

All vector spaces are assumed to be complex throughout this paper. Given a vector space $V$, we let $M_{m, n}(V)$ denote the matrix
space of all $m$ by $n$ matrices $v=[v_{ij}]$ with $v_{ij}\in V$, and we set $M_n(V)=M_{n,n}(V)$. If $V=\mathbb C$,  we write
$M_{m,n}=M_{m,n} (\mathbb C)$ and $M_n=M_{n,n}(\mathbb C)$, which means that we may identify $M_{m,n}(V)$ with the tensor product
$M_{m,n}\otimes V$. We identify $M_{m,n}$ with the normed space $\mathcal B({\mathbb C}^n, {\mathbb C}^m)$. We use the standard
matrix multiplication and *-operation for compatible scalar
matrices, and $1_n$ for the identity matrix in $M_n$, and $0_{m,n}$ for the $m$ by $n$ zero matrix.
There are two natural operations on the matrix spaces. For $v\in M_{m,n}(V)$ and $w\in M_{p,q} (V)$, the direct sum $v\oplus w\in
M_{m+p,n+q}(V)$ is defined by letting
\[v\oplus w=\left[\begin{array}{cc}v&0\\ 0&w\end{array}\right],\]
and if we are given $\alpha\in M_{m,p}$, $v\in M_{p,q}(V)$ and $\beta\in M_{q,n}$, the matrix
product  $\alpha v\beta\in M_{m,n}(V)$ is defined by
\[\alpha v\beta=\left[\sum_{k,l}\alpha_{ik}v_{kl}\beta_{lj}\right].\]

A *-{\it vector space} $V$ is a complex vector space together with a conjugate linear mapping $v\longmapsto v^*$ such that
$v^{**}=v$. A *-vector space $V$ is said to be {\it matrix ordered} if:
\begin{enumerate}
\item each $M_n(V)$, $n\in\mathbb N$, is partially ordered;
\item  $\gamma^*M_n(V)^+\gamma\subseteq M_m(V)^+$ if $\gamma=[\gamma_{ij}]$ is any $n\times m$
matrix of complex numbers.
\end{enumerate}
A {\it matrix order unit space} $(A, 1)$ is a matrix ordered space $A$ together with a distinguished order unit $1$ satisfying the following conditions:
\begin{enumerate}
\item $A^+$ is a proper cone with the order unit $1$;
\item each of the cones $M_n(A)^+$ is Archimedean.
\end{enumerate}
Each matrix order unit space $(A, 1)$ may be provided with the norm
\[\|a\|=\inf\left\{t\in\mathbb R:\, \left[\begin{array}{cc} t1&a\\ a^*&t1\end{array}
\right]\ge 0\right\}.\]
In this paper, we will assume that $A$ is complete for the norm. For a matrix order unit space $(A, 1)$, The {\it
matrix state space} of $(A, 1)$ is the collection $\mathcal{CS}(A)=(CS_n(A))$ of {\it matrix states} $CS_n(A)=\{\varphi: 
\varphi$ is a unital completely positive linear mapping from $A$ into $M_n\}$.

If $V$ and $W$ are vector spaces in duality, then they determine
the matrix pairing
\[\ll\cdot, \cdot\gg:\, M_n(V)\times M_m(W)\mapsto M_{nm},\]
where \[\ll[v_{ij}], [w_{kl}]\gg=\left[<v_{ij}, w_{kl}>\right]\]
for $[v_{ij}]\in M_n(V)$ and $[w_{kl}]\in M_m(W)$.

A {\it graded set} $\mathbf S=(S_n)$ is a sequence of sets $S_n(n\in\mathbb N)$. If, for each $n\in\mathbb{N}$, $M_n(V)$ is a topological 
space, a graded set $\mathbf S=(S_n)$ with $S_n\subseteq M_n(V)$ is {\it closed} or {\it compact} if
that is the case for each set $S_n$ in the topology on $M_n(V)$. Given a vector space $V$, we say that a graded set
$\mathbf B=(B_n)$ with $B_n\subseteq M_n(V)$ is {\it absolutely
matrix convex} if for all $m, n\in\mathbb N$, $B_m\oplus B_n\subseteq B_{m+n}$, and $\alpha B_m\beta\subseteq B_n$ for any contractions 
$\alpha\in M_{n,m}$ and $\beta\in M_{m,n}$. A {\it matrix convex set} in $V$ is a graded set
$\mathbf K=(K_n)$ of subsets $K_n\subseteq M_n(V)$ such that $\sum^k_{i=1}\gamma_i^\ast v_i\gamma_i\in K_n$ for all $v_i\in K_{n_i}$ and 
$\gamma_i\in M_{n_i,n}$ for $i=1,2,\cdots,k$ satisfying $\sum^k_{i=1}\gamma_i^\ast\gamma_i= 1_n$.
Let $V$ and $W$ be vector spaces in duality, and let $\mathbf
S=(S_n)$ be a graded set with $S_n\subseteq M_n(V)$. The {\it
absolute operator polar} ${\mathbf
S}^{\circledcirc}=(S_n^{\circledcirc})$ with
$S_n^\circledcirc\subseteq M_n(W)$, is defined by
$S_n^\circledcirc=\{w\in M_n(W):\,\|\ll v, w\gg\|\le 1\hbox{ for
all }v\in S_r, r\in\mathbb N\}$. 

Given an arbitrary vector space $V$, a {\it matrix gauge} $\mathcal G=(g_n)$ on $V$ is a sequence of gauges $g_n: M_n(V)\mapsto 
[0, +\infty]$ such that
\begin{enumerate}
\item $g_{m+n}(v\oplus w)=\max\{g_m(v), g_n(w)\}$;
\item $g_n(\alpha v\beta)\le\|\alpha\| g_m(v)\|\beta\|$,
\end{enumerate}
for any $v\in M_m(V)$, $w\in M_n(V)$, $\alpha\in M_{n,m}$ and $\beta\in M_{m,n}$. A matrix gauge $\mathcal G=(g_n)$ is a {\it
matrix seminorm} on $V$ if for any $n\in\mathbb N, g_n(v)<+\infty$ for all $v\in M_n(V)$. If each $g_n$ is a norm on $M_n(V)$, we say
that $\mathcal{G}$ is a {\it matrix norm}. An {\it operator space} is a vector space together with a matrix norm on it. Given two operator 
spaces $V$ and $W$. We denote by $\mathrm{CB}(V,W)$ the Banach space of all completely bounded linear mapping $\varphi$ from $V$ into $W$ 
equipped with the completely bounded norm $\|\varphi\|_{cb}$. An {\it operator system} is a 
closed unital self-adjoint linear subspace of a unital $C^\ast$-algebra. For a matrix
order unit space $(A, 1)$, it is an operator space with the matrix norm determined by the matrix order on it. Every matrix order unit space is
completely order isomorphic to an operator system\cite{chef}.

\section{Continuous mappings on matrix metric spaces}\label{se:3}

First we recall the definition of matrix metrics on graded sets. See \cite{wu2} for more details.

\begin{definition}\label{de:31} 
Let $V$ be a vector space and let $\mathbf{K}=(K_n)$ be a graded set with
$K_n\subseteq M_n(V)$. A {\it matrix metric} ${\mathcal D}=(D_n)$ on $\mathbf K$ is a sequence
of metrics
\[D_n:\, K_n\times K_n\longmapsto [0, +\infty)\]
such that
\begin{enumerate}
\item if $x, u\in K_m$ and $y, v\in K_n$ such that $x\oplus y, u\oplus v\in K_{m+n}$, then
$D_{m+n}(x\oplus y, u\oplus v)=\max\{D_m(x, u), D_n(y, v)\}$;
\item if $x, u\in K_m$ and $\alpha\in M_{m, n}$ with $\alpha^*\alpha=1_n$ such that
$\alpha^* x\alpha, \alpha^* u\alpha\in K_n$, then $D_n(\alpha^*
x\alpha, \alpha^* u\alpha) \le D_m(x, u)$.
\end{enumerate}
The ordered pair $(\mathbf{K},{\mathcal D})$ is said to be a {\it matrix metric space} over $V$.
\end{definition}

\begin{example}\label{ex:32}
Let $V$ be an operator space with matrix norm $\|\cdot\|=(\|\cdot\|_n)$ and let $\mathbf{K}=(K_n)$ be a graded set with
$K_n\subseteq M_n(V)$. For $n\in\mathbb{N}$ and $x,y\in K_n$, we define 
\[D_n(x,y)=\|x-y\|_n.\]
Then $\mathcal{D}=(D_n)$ is a matrix metric on $\mathbf{K}$, and is called the matrix metric induced by the matrix norm on $V$. 
$(\mathbf{K},\mathcal{D})$ is called a matrix metric space over the operator space $V$.
\end{example}

\begin{example}\label{ex:33}
Assume $\mathcal L=(L_n)$ is a matrix Lipschitz seminorm on the
matrix order unit space $(\mathcal V, 1)$ (here we do not assume that $\mathcal{V}$ is complete for the matrix norm $\|\cdot\|=(\|\cdot\|_n)$ 
determined by the matrix order on it) and the image of $L_1^1=\{a\in\mathcal{V}: L_1(a)\le 1\}$ in 
$\tilde{\mathcal{V}}=\mathcal{V}/(\mathbb{C}1)$ is totally
bounded for $\|\cdot\|_1^\sim$. Then the sequence $\mathcal{D}_{\mathcal{L}}=(D_{L_n})$ of
metrics defined by
\[D_{L_n}(\varphi,\psi)=\sup\{\|\ll\varphi,a\gg-\ll\psi,a\gg\|: L_r(a)\le 1, r\in\mathbb{N}\},\]
for $\varphi,\psi\in CS_n(\mathcal{V})$ and $n\in\mathbb{N}$, is a matrix
metric on $\mathcal{CS(V)}$(see Theorem 5.3 in \cite{wu1}).
\end{example}

The theory of operator spaces is closely related to the structure of matrices over the spaces and the mappings on them. Here we extend 
the concept of mappings to this situation.

\begin{definition}\label{de:34} 
Let $V$ and $W$ be two vector spaces and let $\mathbf{K}=(K_n)$ and $\mathbf{G}=(G_n)$ be two graded sets with
$K_n\subseteq M_n(V)$ and $G_n\subseteq M_n(G)$. A {\it matrix mapping} from $\mathbf{K}$ into $\mathbf{G}$ is a sequence 
$\mathbf{f}=(f_n)$ of mappings $f_n: K_n\mapsto G_n$. We denote by $\mathcal{M}(\mathbf{K},\mathbf{G})$ the set of all matrix mappings 
from $\mathbf{K}$ into $\mathbf{G}$.
\end{definition}

When each $G_n$ is a subspace of $M_n(G)$, we can define
\[\mathbf{f}+\mathbf{g}=(f_n+g_n),\ \ \ \alpha\mathbf{f}=(\alpha f_n),\]
for $\mathbf{f}=(f_n), \mathbf{g}=(g_n)\in\mathcal{M}(\mathbf{K},\mathbf{G})$ and $\alpha\in\mathbb{C}$. Then 
$\mathcal{M}(\mathbf{K},\mathbf{G})$ is a vector space over $\mathbb{C}$. 

Now we define the matrix analogue of the equicontinuity.

\begin{definition}\label{de:35} 
Given two matrix metric spaces $(\mathbf{K}_1,\mathcal{D}_1)$ and $(\mathbf{K}_2,\mathcal{D}_2)$, we say that a matrix mapping 
$\mathbf{f}=(f_n)$ from $\mathbf{K}_1$ into $\mathbf{K}_2$ is {\it continuous} if each $f_n$ is continuous. 
We let $\mathcal{C}(\mathbf{K}_1,\mathbf{K}_2)$ denote the space of all continuous mappings 
$\mathbf{f}: \mathbf{K}_1\mapsto \mathbf{K}_2$. A subset $S\subseteq\mathcal{C}(\mathbf{K}_1,\mathbf{K}_2)$ is said to be 
{\it equicontinuous} if for 
any $\epsilon>0$, $n\in\mathbb{N}$ and $v\in K_{1,n}$ there exists a $\delta=\delta(\epsilon,v,n)>0$ such that 
\[D_{2,n}(f_n(v),f_n(w))<\epsilon,\]
for all $\mathbf{f}=(f_k)\in S$ and $w\in K_{1,n}$ with $D_{1,n}(v,w)<\delta$. A subset $S\subseteq\mathcal{C}(\mathbf{K}_1,\mathbf{K}_2)$ is 
said to be {\it uniformly equicontinuous} if for 
any $\epsilon>0$ there exists a $\delta=\delta(\epsilon)>0$ such that 
\[D_{2,n}(f_n(v),f_n(w))<\epsilon,\]
for all $\mathbf{f}=(f_k)\in S$ and $n\in\mathbb{N}$ and $w,v\in K_{1,n}$ with $D_{1,n}(v,w)<\delta$.
\end{definition}

A matrix mapping $\mathbf{f}=(f_n)$ from matrix convex set $\mathbf{K}=(K_n)$ into matrix convex set $\mathbf{G}=(G_n)$ is said to be 
{\it matrix affine} if $f_n(\sum^k_{i=1}\gamma_i^\ast v_i\gamma_i)=\sum^k_{i=1}\gamma_i^\ast f_{n_i}(v_i)\gamma_i$ for all $v_i\in K_{n_i}$ 
and $\gamma_i\in M_{n_i,n}$ for $i=1,2,\cdots,k$ satisfying $\sum^k_{i=1}\gamma_i^\ast\gamma_i= 1_n$\cite{wewi}. We let $A(\mathbf{K}, 
\mathbf{G})$ denote the set of all matrix affine mappings from $\mathbf{K}$ into $\mathbf{G}$. Clearly the uniform equicontinuity of a 
subset $S\subseteq\mathcal{C}(\mathbf{K}_1,\mathbf{K}_2)$ implies the equicontinuity of it. The following proposition indicates its 
converse holds in some special cases.

\begin{proposition}\label{pr:36}
Suppose that $(\mathbf{K},\mathcal{D})$ is a matrix metric space and $\mathbf{K}$ is compact matrix convex. Let $(\mathbf{G},\mathcal{E}^{(k)})$ 
be the matrix metric space with $G_n=M_n(M_k)$ and the matrix metric $\mathcal{E}^{(k)}=(E^{(k)}_n)$ induced by the matrix norm on $M_k$. Denote by 
$A\mathcal{C}(\mathbf{K},\mathbf{G})$ the set of all continuous and matrix affine mappings from $\mathbf{K}$ into $\mathbf{G}$. If 
$S\subseteq A\mathcal{C}(\mathbf{K},\mathbf{G})$ is equicontinuous, the $S$ is uniformly equicontinuous.
\end{proposition}

\begin{proof}
Given $\epsilon>0$. Since $S$ is equicontinuous, for any $v\in K_k$ we can find $\delta=\delta(\epsilon,v)>0$ such that 
\[E_k^{(k)}(f_k(v),f_k(w))<\frac{\epsilon}8,\]
for all $\mathbf{f}=(f_n)\in S$ and $w\in K_k$ with $D_k(v,w)<\delta$. That $K_k$ is compact means that there are $v_1,\cdots,v_t\in K_k$ 
such that $K_k\subseteq\cup_{i=1}^tU(v_i;\frac{\delta_i}2)$, where $U(v_i;\frac{\delta_i}2)=\{w\in K_k:D_k(v_i,w)<\frac{\delta_i}2\}$ and 
$\delta_i=\delta(\epsilon, v_i)$. Take $\delta_0=\frac12\min\{\delta_1,\cdots,\delta_t\}$. Then for $w,v\in K_k$ with $D_k(v,w)<\delta_0$, 
there is a $v_i$ with $v\in U(v_i;\frac{\delta_i}2)$. So $D_k(w,v_i)\le D_k(w,v)+D_k(v,v_i)<\delta_i$. Hence for any $\mathbf{f}=(f_n)\in S$ 
we have
\[E_k^{(k)}(f_k(v),f_k(w))\le E_k^{(k)}(f_k(v),f_k(v_i))+E_k^{(k)}(f_k(v_i),f_k(w))<\frac{\epsilon}4.\]

Suppose that $v,w\in K_r$ with $D_r(v,w)<\delta_0$ and $r\neq k$. If $r<k$, choose $u\in K_{k-r}$. Then 
\[v\oplus u=[1_r\ 0_{r,k-r}]^\ast v [1_r\ 0_{r,k-r}]+[0_{k-r,r}\ 1_{k-r}]^\ast u[0_{k-r,r}\ 1_{k-r}]\in K_k\]
since $[1_r\ 0_{r,k-r}]^\ast [1_r\ 0_{r,k-r}]+[0_{k-r,r}\ 1_{k-r}]^\ast [0_{k-r,r}\ 1_{k-r}]=1_k$ and $\mathbf{K}$ is matrix convex. Similarly, $w\oplus u\in K_k$. 
Also $D_k(v\oplus u,w\oplus u)=D_r(v,w)<\delta_0$. So for any $\mathbf{f}=(f_n)\in S$ we have 
\[E_k^{(k)}(f_k(v\oplus u), f_k(w\oplus u))<\frac{\epsilon}4.\]
But that $\mathbf{f}$ is matrix affine implies
\[f_k(v\oplus u)=f_r(v)\oplus f_{k-r}(u),\ \ \ f_k(w\oplus u)=f_r(w)\oplus f_{k-r}(u),\]
that is, $E_r^{(k)}(f_r(v),f_r(w))=E_k^{(k)}(f_r(v)\oplus f_{k-r}(u),f_r(w)\oplus f_{k-r}(u))=E_k^{(k)}(f_k(v\oplus u), f_k(w\oplus u))<\frac{\epsilon}4<
\epsilon$. If $r>k$, for any unit vector $\xi\in\mathbb{C}^r\otimes\mathbb{C}^k$ there exist an isometry $\alpha: \mathbb{C}^k\mapsto\mathbb{C}^r$ and a unit vector 
$\xi_1\in\mathbb{C}^k\otimes\mathbb{C}^k$ such that $\xi=(\alpha\otimes 1_k)(\xi_1)$ by Lemma 5.1 in \cite{efwe}. Since $\mathbf{K}$ is matrix convex and 
$D_k(\alpha^\ast v\alpha,\alpha^\ast w\alpha)\le D_r(v,w)<\delta_0$, we have
\[\begin{array}{rcl}
&&|<(f_r(v)-f_r(w))\xi,\xi>|\\
&=&|<(f_r(v)-f_r(w))(\alpha\otimes 1_k)(\xi_1),(\alpha\otimes 1_k)(\xi_1)>|\\
&=&|<(\alpha^\ast\otimes 1_k)(f_r(v)-f_r(w))(\alpha\otimes 1_k)(\xi_1),\xi_1>|\\
&=&|<(f_k(\alpha^\ast v\alpha)-f_k(\alpha^\ast w\alpha))\xi_1,\xi_1>|\\
&\le&\|f_k(\alpha^\ast v\alpha)-f_k(\alpha^\ast w\alpha)\|\\
&=&E_k^{(k)}(f_k(\alpha^\ast v\alpha),f_k(\alpha^\ast w\alpha))\\
&<&\frac{\epsilon}4.
\end{array}\]
For any unit vectors $\eta,\zeta\in\mathbb{C}^r\otimes\mathbb{C}^k$, we have
\[\begin{array}{rcl}
&&|<(f_r(v)-f_r(w))\eta,\zeta>|\\
&\le&\left|<(f_r(v)-f_r(w))\frac{\eta+\zeta}2,\frac{\eta+\zeta}2>\right|\\
&&+\left|<(f_r(v)-f_r(w))\frac{\eta-\zeta}2,\frac{\eta-\zeta}2>\right|\\
&&+\left|i<(f_r(v)-f_r(w))\frac{\eta+i\zeta}2,\frac{\eta+i\zeta}2>\right|\\
&&+\left|i<(f_r(v)-f_r(w))\frac{\eta-i\zeta}2,\frac{\eta-i\zeta}2>\right|\\
&<&\epsilon.
\end{array}\]
From the arbitrariness of $\eta$ and $\zeta$, we obtain that $E_r^{(k)}(f_r(v),f_r(w))=\|f_r(v)-f_r(w)\|\le\epsilon$. By definition, $S$ is 
uniformly equicontinuous.
\end{proof}

\section{Matrix metric sets}\label{se:4}

Motivated by the idea of Antonescu and Christensen and our results in \cite{wu2}, we give the operator space version of the metric set.
 
\begin{definition}\label{de:41}
Let $(A,1)$ be a matrix order unit space. A graded set $\mathbf{K}=(K_n)$ with $K_n\subseteq M_n(A)$ is called a 
{\it matrix metric set} of $(A,1)$ if it is norm compact, self-adjoint and absolutely matix convex, and separates the matrix states on $(A,1)$.
\end{definition}

One can easily construct matrix metric sets for separable matrix order unit spaces.

\begin{example}\label{ex:415}
Given a countable group $G=\{g_n: n\in\mathbb{N}\}$ and a closed self-adjoint subspace $A$ of $C^\ast_r(G)$ containing the unit 
$\lambda_e$, where $e$ is the identity element of $G$. Then with the usual partial ordering on $M_n(A)$ for $n\in\mathbb{N}$, 
$(A, \lambda_e)$ becomes a matrix order unit space. Set
\[K=A\cap\overline{\mathrm{co}}\left(\cup_{n=1}^\infty\{\alpha\lambda_{g_n}+\beta\lambda^\ast_{g_n}: |\alpha|+|\beta|\le\frac1n, 
\alpha,\beta\in\mathbb{C}\}\right),\]
where $\overline{\mathrm{co}}$ means the closed convex hull. Then $K$ is a norm compact, self-adjoint and absolutely convex subset of $A$, and so is weakly 
closed and absolutely convex. Thus there is a weakly closed absolutely matrix convex set $\mathbf{K}=(K_n)$ with $K_n\subseteq M_n(A)$ 
and $K_1=K$ (see page 181 in \cite{efwe}). Clearly $\mathbf{K}$ is norm closed. For any $n\in\mathbb{N}$ and $a=[a_{ij}]\in K_n$, we have
\[a_{ij}=[0\ \cdots \ 0 \ 1_i \ 0 \cdots 0] a [0\ \cdots \ 0 \ 1_j \ 0 \cdots 0]^\ast\in K_1.\]
Since $K_1=K$ is norm compact, $K_1$ is totally bounded, and hence every $K_n$ is totally bounded. Each $K_n$ is norm closed implies 
it is also norm compact. Clearly $K$ separates the states on $A$. Since the matrix state space $\mathcal{CS}(A)$ of $A$ is matrix convex, $K$ 
also separates the matrix states on $A$. So $\mathbf{K}$ separates the matrix states on $A$. Therefore, $\mathbf{K}$ is a matrix metric set 
of $(A, \lambda_e)$.
\end{example}

The natural topology on the matrix state space $\mathcal{CS}(A)$ is the BW-topology, that is, topologies each $CS_n(A)$ by BW-topology (see 
page 146 in \cite{ar}). The following result justifies the definition of a matrix metric set, that is, it generates the BW-topology.

\begin{proposition}\label{pr:42}
Let $(A,1)$ be a matrix order unit space, $\mathcal{CS}(A)$ the matrix state space of $(A,1)$ and $\mathbf{K}=(K_n)$ a matrix 
metric set of $(A,1)$. Then $\mathcal{D}_{\mathcal{K}}=(D_{K_n})$, where 
\[D_{K_n}(\varphi,\psi)=\sup\{\|\ll\varphi,a\gg-\ll\psi,a\gg\|: a\in K_r, r\in\mathbb{N}\},\]
for $\varphi,\psi\in CS_n(A)$ and $n\in\mathbb{N}$, is a matrix metric on $\mathcal{CS}(A)$ and the $\mathcal{D}_{\mathcal{K}}$-
topology on $\mathcal{CS}(A)$ agrees with the BW-topology.
\end{proposition}

\begin{proof}
Denote
\[C_n(\varphi,\psi)=\sup\{\|\ll\varphi,a\gg-\ll\psi,a\gg\|: a\in K_n\},\]
for $\varphi,\psi\in CS_n(A)$ and $n\in\mathbb{N}$. Clearly $C_n(\varphi,\psi)\le D_{K_n}(\varphi,\psi)$ for all $\varphi,\psi\in CS_n(A)$ 
and $n\in\mathbb{N}$. 

For $a=[a_{pq}]\in K_r$ and $r<n$, we have that $0_{n-r}=0_{n-r,r}a0_{r,n-r}\in K_{n-r}$ and hence $a\oplus 0_{n-r}\in K_n$ because $\mathbf{K}$ 
is absolutely matrix convex. So for $\varphi,\psi\in CS_n(A)$, we have 
\[\begin{array}{rcl}
&&\|\ll\varphi,a\gg-\ll\psi,a\gg\|\\
&=&\|\ll\varphi,a\oplus 0_{n-r}\gg-\ll\psi,a\oplus 0_{n-r}\gg\|\le C_n(\varphi,\psi).
\end{array}\]
Suppose that $r>n$. For arbitrary unit vectors $\xi,\eta\in\mathbb{C}^r\otimes\mathbb{C}^n$, there exist isometries $\alpha,\beta:\mathbb{C}^n
\mapsto\mathbb{C}^r$ and unit vectors $\xi_1,\eta_1\in\mathbb{C}^n\otimes\mathbb{C}^n$ for which $\xi=(\alpha\otimes 1_n)(\xi_1)$ and 
$\eta=(\beta\otimes 1_n)(\eta_1)$ by Lemma 5.1 in \cite{efwe}. That $\mathbf{K}$ 
is absolutely matrix convex implies
\[\begin{array}{rcl}
&&|<(\ll\varphi,a\gg-\ll\psi,a\gg)\eta,\xi>|\\
&=&|<(\ll\varphi,a\gg-\ll\psi,a\gg)(\beta\otimes 1_n)(\eta_1),(\alpha\otimes 1_n)(\xi_1)>|\\
&=&|<(\ll\varphi,\alpha^\ast a\beta\gg-\ll\psi,\alpha^\ast a\beta\gg)\eta_1,\xi_1>|\\
&\le&\|\ll\varphi,\alpha^\ast a\beta\gg-\ll\psi,\alpha^\ast a\beta\gg\|\\
&\le&C_n(\varphi,\psi).
\end{array}\]
Since $\xi$ and $\eta$ are arbitrary unit vectors, we conclude that $\|\ll\varphi,a\gg-\ll\psi,a\gg\|\le C_n(\varphi,\psi)$. Therefore, 
$D_{K_n}(\varphi,\psi)\le C_n(\varphi,\psi)$ for all $\varphi,\psi\in CS_n(A)$ and $n\in\mathbb{N}$, and so $C_n(\varphi,\psi)=
D_{K_n}(\varphi,\psi)$ for all $\varphi,\psi\in CS_n(A)$ and $n\in\mathbb{N}$. Since $\mathbf{K}$ separates $\mathcal{CS}(A)$ and norm compact, each $D_{K_n}$ is a bounded metric on $CS_n(A)$.

For $\varphi_1,\varphi_2\in CS_m(A)$ and $\psi_1,\psi_2\in CS_p(A)$, we have that $\varphi_1\oplus\psi_1,\varphi_2\oplus\psi_2\in 
CS_{m+p}(A)$ and
\[\begin{array}{rcl}
&&D_{K_{m+p}}(\varphi_1\oplus\psi_1,\varphi_2\oplus\psi_2)\\
&=&\sup\{\|\ll\varphi_1\oplus\psi_1,a\gg-\ll\varphi_2\oplus\psi_2,a\gg\|: a\in K_r, r\in\mathbb{N}\}\\
&=&\sup\{\max\{\|\ll\varphi_1,a\gg-\ll\varphi_2,a\gg\|, \\
&&\|\ll\psi_1,a\gg-\ll\psi_2,a\gg\|\}: a\in K_r, r\in\mathbb{N}\}\\
&=&\max\{\sup\{\|\ll\varphi_1,a\gg-\ll\varphi_2,a\gg\|, a\in K_r, r\in\mathbb{N}\}, \\
&&\sup\{\|\ll\psi_1,a\gg-\ll\psi_2,a\gg\|: a\in K_r, r\in\mathbb{N}\}\}\\
&=&\max\{D_{K_m}(\varphi_1,\varphi_2), D_{K_p}(\psi_1,\psi_2)\}.
\end{array}\]
If $\varphi,\psi\in CS_m(A)$ and $\alpha\in M_{m,p}$ with $\alpha^\ast\alpha=1_n$, then $\alpha^\ast\varphi\alpha, 
\alpha^\ast\psi\alpha\in CS_p(A)$, and 
\[\begin{array}{rcl}
&&D_{K_{p}}(\alpha^\ast\varphi\alpha,\alpha^\ast\psi\alpha)\\
&=&\sup\{\|\ll\alpha^\ast\varphi\alpha,a\gg-\ll\alpha^\ast\psi\alpha,a\gg\|: a\in K_r, r\in\mathbb{N}\}\\
&=&\sup\{\|(\alpha^\ast\otimes 1_r)(\ll\varphi,a\gg-\ll\psi,a\gg)(\alpha\otimes 1_r)\|: a\in K_r, r\in\mathbb{N}\}\\
&\le&\sup\{\|\ll\varphi,a\gg-\ll\psi,a\gg\|, a\in K_r, r\in\mathbb{N}\}\\
&=&D_{K_m}(\varphi,\psi)\}.
\end{array}\]
Therefore, $\mathcal{D}_{\mathcal{K}}$ is a matrix metric on $\mathcal{CS}(A)$. 

Clearly, the topology on each $CS_n(A)$ induced by $D_{K_n}$ is a Hausdorff topology. Suppose $\{\varphi_i\}\subseteq CS_n(A)$, $\varphi\in 
CS_n(A)$ and $\lim_i\varphi_i=\varphi$ in the BW-topology. Then $\lim_i\varphi_i(a)=\varphi(a)$ for all $a\in A$. Given $\epsilon>0$. 
For $a\in K_1$, there is an $I_a$ such that $\|\varphi_i(a)-\varphi(a)\|<\frac{\epsilon}{3n^2}$ for $i\ge I_a$. When $b\in U(a;
\frac{\epsilon}{3n^2})=\{c\in A:\|c-a\|<\frac{\epsilon}{3n^2}\}$, we have
\[\begin{array}{rcl}
&&\|\varphi_i(b)-\varphi(b)\|\\
&\le&\|\varphi_i(b)-\varphi_i(a)\|+\|\varphi_i(a)-\varphi(a)\|+\|\varphi(a)-\varphi(b)\|\\
&\le&2\|b-a\|+\|\varphi_i(a)-\varphi(a)\|\\
&<&\frac{\epsilon}{n^2},
\end{array}\]
for $i\ge I_a$. Since $K_1$ is norm compact, there exists an $s\in\mathbb{N}$ such that $K_1\subseteq\cup_{i=1}^sU(a_i;\frac{\epsilon}{3n^2})$ 
for some $a_1,\cdots,a_s\in K_1$. Then for $i\ge I_{a_j}, j=1,2,\cdots,s$, and $c\in K_1$, we can find an $i_0\in\{1,2,\cdots,s\}$ such that 
$c\in U(a_{i_0},\frac{\epsilon}{3n^2})$, and so $\|\varphi_i(c)-\varphi(c)\|<\frac{\epsilon}{n^2}$.

For $a=[a_{pq}]\in K_n$ and $i\ge I_{a_j}, j=1,2,\cdots,s$, 
we obtain
\[\begin{array}{rcl}
\|\ll\varphi_i,a\gg-\ll\varphi,a\gg\|&=&\|[\varphi_i(a_{pq})-\varphi(a_{pq})]\|\\
&\le&\sum_{p,q=1}^n\|\varphi_i(a_{pq})-\varphi(a_{pq})\|<\epsilon,
\end{array}\]
Therefore, $D_{K_n}(\varphi_i,\psi)\le\epsilon$ for $i\ge I_{a_j}, j=1,2,\cdots,s$, that is, $\lim_i\varphi_i=\varphi$ 
in the $D_{K_n}$-topology. So $\mathcal{D}_{\mathcal{K}}$-topology on $\mathcal{CS}(A)$ is weaker than the BW-topology.

On the other hand, $\mathcal{CS}(A)$ is BW-compact by Theorem 6.4 in \cite{pau}, and so $\mathcal{D}_{\mathcal{K}}$-topology and BW-topology agree.
\end{proof}

From the proof of Proposition \ref{pr:42}, we have

\begin{corollary}\label{co:43}
Let $(A,1)$ be a matrix order unit space, $\mathcal{CS}(A)$ the matrix state space of $(A,1)$ and $\mathbf{K}=(K_n)$ a matrix 
metric set of $(A,1)$. Then  
\[D_{K_n}(\varphi,\psi)=\sup\{\|\ll\varphi,a\gg-\ll\psi,a\gg\|: a\in K_n\},\]
for $\varphi,\psi\in CS_n(A)$ and $n\in\mathbb{N}$.
\end{corollary}

By a {\it matrix Lip-gauge} on a matrix order unit space $(A, 1)$ we mean a
matrix gauge $\mathcal{G}=(G_n)$ on $(A,1)$ such that: (1) the null space of each $G_n$ is $M_n(\mathbb{C}1)$; (2) 
$G_n(v^*)=G_n(v)$ for any $v\in M_n(A)$; (3) $\{v\in A: G_1(v)<+\infty\}$ is dense in $A$; (4) the $\mathcal{D}_{\mathcal{G}}$-topology on 
$\mathcal{CS}(A)$ agrees with the BW-topology. The matrix Lip-gauge $\mathcal{G}=(G_n)$ is {\it lower semicontinuous} if each $G_n$ is 
lower semicontinuous.

\begin{corollary}\label{co:44}
Let $(A,1)$ be a matrix order unit space, $\mathcal{CS}(A)$ the matrix state space of $(A,1)$ and $\mathbf{K}=(K_n)$ a matrix 
metric set of $(A,1)$. Then there is a lower semicontinuous matrix Lip-gauge $\mathcal{L}=(L_{n})$ on $(A,1)$ such that
\[D_{K_n}(\varphi,\psi)=\sup\{\|\ll\varphi,a\gg-\ll\psi,a\gg\|: L_r(a)\le 1, r\in\mathbb{N}\},\]
for $\varphi,\psi\in CS_n(A)$ and $n\in\mathbb{N}$.
\end{corollary}

\begin{proof}
Clearly $\mathcal{D}_{\mathcal{K}}$ is convex, midpoint balance, and midpoint concave. Now the corollary follows from Theorem 6.12 in \cite{wu2}.
\end{proof}

\section{Relative compactness}\label{se:5}

In this section we state and prove our main results. First let us take a look at what is the boundedness of matrix mappings.

\begin{definition}\label{de:51}
Given two matrix metric spaces $(\mathbf{K}_1,\mathcal{D}_1)$ and $(\mathbf{K}_2,\mathcal{D}_2)$. Fix an $x=(x_n)\in\mathbf{K}_2$ (that is, 
each $x_n\in K_{2,n}$). For a matrix mapping $\mathbf{f}=(f_n)\in\mathcal{M}(\mathbf{K}_1,\mathbf{K}_2)$, we define
\[p(\mathbf{f})=\sup\{D_{2,n}(f_n(w),x_n): w\in K_{1,n}, n\in\mathbb{N}\}.\]
If $p(\mathbf{f})<+\infty$, we say that $\mathbf{f}$ is {\it bounded}. Given $S\subseteq\mathcal{M}(\mathbf{K}_1,\mathbf{K}_2)$. If there is a 
constant $C>0$ such that $p(\mathbf{f})\le C$ for any $\mathbf{f}\in S$, we say that $S$ is {\it bounded}.
\end{definition}

If $(\mathbf{K}_2,\mathcal{D}_2)$ is the matrix metric space over an operator space, we take $x=(0_n)$. Then clearly 
\[p(\mathbf{f})=\sup\{\|f_n(w)\|_n: w\in K_{1,n}, n\in\mathbb{N}\}\]
is a faithful gauge on $\mathcal{M}(\mathbf{K}_1,\mathbf{K}_2)$. The following lemma displays one of their aspects of the 
boundedness of matrix mappings.

\begin{lemma}\label{le:52}
Let $(A,1)$ be a matrix order unit space. The canonical mapping of $A$ into $A(\mathcal{CS}(A))$, which sends $x\in M_n(A)$ to 
$\hat{x}=(\hat{x}_r)\in M_n(A(\mathcal{CS}(A)))\simeq A(\mathcal{CS}(A),M_n)$ given by $\hat{x}_r(\varphi)=\ll\varphi,x\gg$ for 
$\varphi\in CS_r(A)$ and $r\in\mathbb{N}$, is a unital matrix order preserving bijection between $A$ and $A(\mathcal{CS}(A))$ $($see 
page $314$ in \cite{wewi}$)$. Define 
\[p_n(\hat{x})=\sup\{\|\ll\varphi,x\gg\|: \varphi\in CS_r(A), r\in\mathbb{N}\},\]
for $x\in M_n(A)$ and $n\in\mathbb{N}$. Then each $p_n$ is a norm on $A(\mathcal{CS}(A),M_n)$. Moreover,
\[p_n(\hat{x})=\|x\|_n,\]
for $x\in M_n(A)$ and $n\in\mathbb{N}$.
\end{lemma}

\begin{proof}
Since 
\[p_n(\hat{x})=\sup\{\|\hat{x}_r(\varphi)\|: \varphi\in CS_r(A), r\in\mathbb{N}\},\]
$p_n$ is a faithful gauge. For $x\in M_n(A)$, we have
\[\begin{array}{rcl}
\|x\|_n
&=&\inf\left\{t\in\mathbb{R}: \left[\begin{array}{cc}t1_n&x\\ x^\ast&t1_n\end{array}\right]\ge 0\right\}\\
&=&\inf\left\{t\in\mathbb{R}: \ll\varphi,\left[\begin{array}{cc}t1_n&x\\ x^\ast&t1_n\end{array}\right]\gg\ge 0, 
\varphi\in CS_r(A),r\in\mathbb{N}\right\}\\
&=&\inf\left\{t\in\mathbb{R}: \left[\begin{array}{cc}t1_r\otimes 1_n&\ll\varphi,x\gg\\ \ll\varphi,x\gg^\ast&t1_r\otimes 1_n
\end{array}\right]\ge 0, \varphi\in CS_r(A),r\in\mathbb{N}\right\}\\
&=&\sup\left\{\|\ll\varphi,x\gg\|: \varphi\in CS_r(A),r\in\mathbb{N}\right\}\\
&=&p_n(\hat{x}).
\end{array}\]
\end{proof}

\begin{lemma}\label{le:53}
With notation as in Lemma $\ref{le:52}$, we have 
\[\frac14 p_n(\hat{x})\le q_n(\hat{x})\le p_n(\hat{x}),\]
for $x\in M_n(A)$ and $n\in\mathbb{N}$, where $q_n(\hat{x})=\sup\{\|\ll\varphi,x\gg\|:\varphi\in CS_n(A)\}$.
\end{lemma}

\begin{proof}
Clearly, $q_n(\hat{x})\le p_n(\hat{x})$. Given $\varphi\in CS_r(A)$. If $r<n$, we have
\[\begin{array}{rcl}
\psi&=&[1_r\ 0_{r,n-r}]^\ast\varphi [1_r\ 0_{r,n-r}]\\
&&+[0_{n-r,r}\ 1_{n-r}]^\ast([1\ 0_{1,r-1}]\varphi [1\ 0_{1,r-1}]^\ast\otimes 
1_{n-r})[0_{n-r,r}\ 1_{n-r}]\in CS_n(A),\end{array}\]
and so
\[q_n(\hat{x})\ge\|\ll\psi,x\gg\|\ge\|\ll\varphi,x\gg\|.\]
If $r>n$, for any unit vector $\xi\in\mathbb{C}^r\otimes\mathbb{C}^n$ there exist an isometry $\alpha: \mathbb{C}^n\mapsto\mathbb{C}^r$ 
and a unit vector $\xi_1\in\mathbb{C}^n\otimes\mathbb{C}^n$ for which $(\alpha\otimes 1_n)(\xi_1)=\xi$. So we get
\[\begin{array}{rcl}
|<\ll\varphi,x\gg\xi,\xi>|&=&|<\ll\varphi,x\gg(\alpha\otimes 1_n)(\xi_1),(\alpha\otimes 1_n)(\xi_1)>|\\
&=&|<\ll\alpha^\ast\varphi\alpha,x\gg\xi_1,\xi_1>|\\
&\le&\|\ll\alpha^\ast\varphi\alpha,x\gg\|\le q_n(\hat{x}).
\end{array}\]
Now for any unit vectors $\eta,\zeta\in\mathbb{C}^r\otimes\mathbb{C}^n$, we have
\[\begin{array}{rcl}
&&|<\ll\varphi,x\gg\eta,\zeta>|\\
&\le&\left|<\ll\varphi,x\gg\frac{\eta+\zeta}2,\frac{\eta+\zeta}2>\right|+\left|<\ll\varphi,x\gg\frac{\eta-\zeta}2,\frac{\eta-\zeta}2>\right|\\
&&+\left|i<\ll\varphi,x\gg\frac{\eta+i\zeta}2,\frac{\eta+i\zeta}2>\right|+
\left|i<\ll\varphi,x\gg\frac{\eta-i\zeta}2,\frac{\eta-i\zeta}2>\right|\\
&\le&4 q_n(\hat{x}).
\end{array}\]
Thus $\|\ll\varphi,x\gg\|\le 4q_n(\hat{x})$. Since $\varphi$ is arbitrary, we obtain $\frac14 p_n(\hat{x})\le q_n(\hat{x})$. 
\end{proof}

For a matrix order unit space $(A,1)$ and $\epsilon>0$, we denote $\mathbb{C}\mathbf{1}=(M_n(\mathbb{C}1))$ and 
$\mathbf{A}_{\epsilon}=(A_{n,\epsilon})$ with $A_{n,\epsilon}=\{a\in 
M_n(A): \|a\|_n\le\epsilon\}$ for $n\in\mathbb{N}$, where $\|\cdot\|=(\|\cdot\|_n)$ is the matrix norm determined by the matrix order 
on $(A,1)$. Given two graded sets $\mathbf{U}=(U_n)$ and $\mathbf{V}=(V_n)$ with $U_n, V_n\in M_n(A)$, and a sequence $\boldsymbol{\lambda}=(\lambda_n)$ with 
$\lambda_n\in\mathbb{C}$, we define
\[\mathbf{U}+\mathbf{V}=(U_n+V_n),\ \ \ \ \boldsymbol{\lambda}\mathbf{U}=(\lambda_n U_n).\]
If all $\lambda_n=\lambda$, we define $\lambda\mathbf{U}=(\lambda U_n)$. Let $(\mathbf{K},\mathcal{D})$ be a matrix metric space and $\mathbf{G}$ 
a graded set with $\mathbf{G}\subseteq\mathbf{K}$, that is, each $G_n\subseteq K_n$. If every $G_n$ is bounded with respect to the metric 
$D_n$, we say that $\mathbf{G}$ is bounded.

\begin{theorem}\label{th:54}
Let $(A,1)$ be a matrix order unit space and $\mathbf{K}$ a matrix metric set of $(A,1)$. For any graded set $\mathbf{S}=(S_n)$ with 
$S_n\subseteq M_n(A)$, the following conditions are 
equivalent:
\begin{enumerate}
\item The graded set $\mathbf{S}$ is norm relatively compact.
\item The graded set $\hat{\mathbf{S}}=(\hat{S}_n)$, where $\hat{S}_n=\{\hat{a}\in A(\mathcal{CS}(A),M_n): a\in S_n\}$, is bounded and 
equicontinuous with respect to the matrix metric $\mathcal{D}_{\mathcal{K}}$ and the matrix metric $\mathcal{E}^{(n)}$ induced by the matrix norm 
on each $M_n$.
\item The graded set $\mathbf{S}$ is bounded and for every $\epsilon>0$ there exists a sequence $\boldsymbol{\lambda}=(\lambda_n)$ with $\lambda_n>0$ such that 
\[\mathbf{S}\subseteq \mathbf{A}_{\epsilon}+\boldsymbol{\lambda}\mathbf{K}+\mathbb{C}\mathbf{1}.\]
\end{enumerate}
\end{theorem}

\begin{proof} We will show $(1)\Longrightarrow (2)\Longrightarrow (3)\Longrightarrow (1)$. 
Assume (1). Then each $S_n$ is norm relatively compact. Since $M_n(A)$ is complete, each $S_n$ is totally bounded. By Lemma \ref{le:52}, each 
$\hat{S}_n$ is totally bounded, and so each $\hat{S}_n$ is bounded.

Given $a=[a_{ij}]\in S_n$ and $\epsilon>0$. For any $k\in\mathbb{N}$, suppose that $\varphi\in CS_k(A), \{\varphi_s\}\subseteq CS_k(A)$ and 
$\lim_s\varphi_s=\varphi$ in the $D_{K_k}$-topology. By Proposition \ref{pr:42}, $D_{K_k}$-topology on $CS_k(A)$ agrees with the BW-topology. 
So there exists an $s_0$ such that when $s>s_0$, we have
\[\|\varphi_s(a_{ij})-\varphi (a_{ij})\|<\frac{\epsilon}{n^2},\]
for $i,j=1,2,\cdots,n$. Now for $s>s_0$, we have
\[\begin{array}{rcl}
E^{(n)}_k(\hat{a}_k(\varphi_s),\hat{a}_k(\varphi))
&=&\|\hat{a}_k(\varphi_s)-\hat{a}_k(\varphi)\|\\
&=&\|\ll\varphi_s,a\gg-\ll\varphi,a\gg\|\\
&=&\|[\varphi_s(a_{ij})-\varphi(a_{ij})]\|\\
&\le&\sum_{i,j=1}^n\|\varphi_s(a_{ij})-\varphi(a_{ij})\|<\epsilon.
\end{array}\]
Hence $\hat{a}\in\mathcal{C}(\mathcal{CS}(A),M_n)$.

Fix $\epsilon>0, k\in\mathbb{N}$ and $\varphi\in CS_k(A)$. Since $\hat{S}_n$ is totally bounded, there is an $\frac{\epsilon}3$-net 
$\{\hat{a}^{(1)},\cdots,\hat{a}^{(t)}\}$ in $\hat{S}_n$. 
$\{\hat{a}^{(1)},\cdots,\hat{a}^{(t)}\}\subseteq\mathcal{C}(\mathcal{CS}(A),M_n)$ implies that there is a $\delta=\delta(\epsilon,\varphi,
k)>0$ such that
\[E^{(n)}_k\left(\hat{a}^{(s)}_k(\varphi),\hat{a}^{(s)}_k(\psi)\right)<\frac{\epsilon}3, \ \ \ \ s=1,2,\cdots,t,\]
for $\psi\in CS_k(A)$ with $D_{K_k}(\varphi,\psi)<\delta$. For any $\hat{a}\in\hat{S}_n$, there is a $\hat{a}^{(s)}\in
\{\hat{a}^{(1)},\cdots,\hat{a}^{(t)}\}$ such that $p_n(\hat{a}-\hat{a}^{(s)})<\frac{\epsilon}3$. So for $\psi\in CS_k(A)$ with 
$D_{K_k}(\varphi,\psi)<\delta$, we have
\[\begin{array}{rcl}
E^{(n)}_k(\hat{a}_k(\varphi),\hat{a}_k(\psi))
&=&\|\hat{a}_k(\varphi)-\hat{a}_k(\psi)\|\\
&\le&\left\|\hat{a}_k(\varphi)-\hat{a}_k^{(s)}(\varphi)\right\|+\left\|\hat{a}_k^{(s)}(\varphi)-\hat{a}_k^{(s)}(\psi)\right\|\\
&&+\left\|\hat{a}_k^{(s)}(\psi)-\hat{a}_k(\psi)\right\|\\
&\le&2\left\|a-a^{(s)}\right\|_n+E^{(n)}_k\left(\hat{a}^{(s)}_k(\varphi),\hat{a}^{(s)}_k(\psi)\right)\\
&=&2p_n\left(\hat{a}-\hat{a}^{(s)}\right)+E^{(n)}_k\left(\hat{a}^{(s)}_k(\varphi),\hat{a}^{(s)}_k(\psi)\right)\\
&<&\epsilon.
\end{array}\]
By definition, $\hat{S}_n$ is equicontinuous, whence (2).

Assume (2). From the boundedness of $\hat{\mathbf{S}}$ it follows that $\mathbf{S}$ is bounded by Lemma \ref{le:52}. Since the 
$\mathcal{D}_{\mathcal{K}}$-topology on 
$\mathcal{CS}(A)$ agrees with the BW-topology (see Proposition \ref{pr:42}), $\mathcal{CS}(A)$ is compact in the $\mathcal{D}_{\mathcal{K}}$-topology. 
Moreover, $\mathcal{CS}(A)$ is matrix convex. Clearly $\hat{S}_n\subseteq A\mathcal{C}(\mathcal{CS}(A), M_n)$. By Proposition \ref{pr:36}, each $\hat{S}_n$ 
is uniformly equicontinuous because each $\hat{S}_n$ is equicontinuous. Given $\epsilon>0$ and $n\in\mathbb{N}$. We can find a $\delta_n>0$ such that
\[E^{(n)}_k(\hat{a}_k(\varphi), \hat{a}_k(\psi))<\epsilon,\]
for $a\in S_n$ and $\varphi,\psi\in CS_k(A)$ with $D_{K_k}(\varphi,\psi)\le\delta_n$ and $k\in\mathbb{N}$. Because $\varphi-\psi\in K_k^{\circledcirc}$ if 
and only if $D_{K_k}(\varphi,\psi)\le 1$, $D_{K_k}(\varphi,\psi)\le\delta_n$ if and only if $\varphi-\psi\in\delta_n(K_k^{\circledcirc})$. So we have 
\[\|\ll\varphi-\psi, a\gg\|<\epsilon,\]
for $a\in S_n$ and $\varphi,\psi\in CS_k(A)$ with $\varphi-\psi\in\delta_n(K_k^{\circledcirc})$ and $k\in\mathbb{N}$. 

For $f\in M_k(A^\ast)\simeq\mathrm{CB}(A,M_k)$ 
with $\|f\|_{cb}\le2$ and $\ll f,c\gg=0_{rk}$ for $c\in M_r(\mathbb{C}1)$ and $r\in\mathbb{N}$, there are $\varphi_1,\varphi_2,\varphi_3,\varphi_4\in CS_k(A)$ 
such that
\[f=\varphi_1-\varphi_2+i(\varphi_3-\varphi_4)\]
by Lemma 4.1 in \cite{wu2}. If $f\in\delta_n(K_k^{\circledcirc})$, then  
\[\begin{array}{rcl}
\|\ll\varphi_1-\varphi_2,a\gg\|&=&\left\|\ll\frac12(f+f^\ast),a\gg\right\|\\
&\le&\frac12(\|\ll f,a\gg\|+\|\ll f^\ast,a\gg\|)\\
&=&\frac12(\|\ll f,a\gg\|+\|\ll f,a^\ast\gg^\ast\|)\\
&\le&\delta_n,\end{array}\]
for $a\in K_r$ and $r\in\mathbb{N}$. We obtain that $\varphi_1-\varphi_2\in\delta_n(K_k^{\circledcirc})$. 
Similarly, we have $\varphi_3-\varphi_4\in\delta_n(K_k^{\circledcirc})$. So 
for $a\in S_n$ we have
\[\|\ll f,a\gg\|<2\epsilon.\]
Denote $\{\mathbb{C}\mathbf{1}\}_k^{\bot}=\{f\in M_k(A^\ast):\ll f,c\gg=0_{kr}$ for $c\in M_r(\mathbb{C}1),r\in\mathbb{N}\}$ and $T_k=(A^\ast)_{k,2}\cap
\{\mathbb{C}\mathbf{1}\}_k^{\bot}\cap\delta_n(K_k^{\circledcirc})$ for $k\in\mathbb{N}$. Let $\mathbf{T}=(T_k)$. We have
\[\begin{array}{rcl}
T_n^{\circledcirc}&=&\{b\in M_n(A):\|\ll f,b\gg\|\le 1, f\in T_k,k\in\mathbb{N}\}\\
&=&\frac1{2\epsilon}\{b\in M_n(A):\|\ll f,b\gg\|\le 2\epsilon, f\in T_k,k\in\mathbb{N}\}.
\end{array}\]
Thus if $a\in S_n$, then $a\in 2\epsilon T_n^{\circledcirc}$. So $S_n\subseteq 2\epsilon T_n^{\circledcirc}$. Set
\[W_k=A_{k,\frac12}\cup M_k(\mathbb{C}1)\cup\frac1{\delta_n}K_k,\ \ \ k\in\mathbb{N},\]
and $\mathbf{W}=(W_k)$. We have
\[W_n^{\circledcirc}=\{f\in M_n(A^\ast):\|\ll f,b\gg\|\le1,b\in W_k,k\in\mathbb{N}\}=T_n.\]
So $S_n\subseteq 2\epsilon W_n^{\circledcirc\circledcirc}$. Since $A_{k,\frac12}$ and $M_k(\mathbb{C}1)$ are norm closed and $\frac1{\delta_n}K_k$ is norm compact, 
$X_k=A_{k,\frac12}+M_k(\mathbb{C}1)+\frac1{\delta_n}K_k$ is norm closed. Clearly $\mathbf{X}=(X_k)$ is absolutely matrix convex. So $\mathbf{X}$ is weakly 
closed, $\mathbf{W}\subseteq\mathbf{X}$ and $\mathbf{X}$ is contained in any absolutely matrix convex set containing $\mathbf{W}$. But the generalized bipolar theorem says that 
$\mathbf{W}^{\circledcirc\circledcirc}$ equals $\overline{\mathrm{amco}}(\mathbf{W})$, the smallest weakly closed absolutely matrix convex set containing $\mathbf{W}$ (see Proposition 4.1 in \cite{efwe}). 
Therefore, $S_n\subseteq2\epsilon W_n^{\circledcirc\circledcirc}=2\epsilon X_n=A_{n,\epsilon}+M_n(\mathbb{C}1)+\lambda_nK_n$, where $\lambda_n=
\frac{2\epsilon}{\delta_n}$, whence (3).

Assume (3). Then $\mathbf{S}$ is bounded, and so $\hat{\mathbf{S}}$ is bounded by Lemma \ref{le:52}. Given $\epsilon>0$. We can find a sequence $\boldsymbol{\lambda}=
(\lambda_n)$ with $\lambda_n>0$ such that
\[\mathbf{S}\subseteq \mathbf{A}_{\frac{\epsilon}{128}}+\boldsymbol{\lambda}\mathbf{K}+\mathbb{C}\mathbf{1}.\]
For $n\in\mathbb{N}$ and $a\in S_n$, there are $b\in A_{n,\frac{\epsilon}{128}}, c\in K_n$ and $\alpha=[\alpha_{ij}]\in M_n$ such 
that $a=b+\lambda_n c+[\alpha_{ij}1]$. Then for $k\in\mathbb{N}$ and $\varphi,\psi\in CS_k(A)$ with $D_{K_k}(\varphi,\psi)<
\frac{\epsilon}{64\lambda_n}$, we have
\[\begin{array}{rcl}
&&E^{(n)}_k(\hat{a}_k(\varphi),\hat{a}_k(\psi))\\
&=&\|\ll\varphi,a\gg-\ll\psi,a\gg\|\\
&\le&\|\ll\varphi-\psi,b\gg\|+|\lambda_n|\|\ll\varphi,c\gg-\ll\psi,c\gg\|\\
&\le&2\|b\|+|\lambda_n|D_{K_k}(\varphi,\psi)<\frac{\epsilon}{32}.
\end{array}\]
$CS_n(A)$ is BW-compact and $CS_n(A)=\cup_{\varphi\in CS_n(A)}U(\varphi;\frac{\epsilon}{64\lambda_n})$. Hence there are $\varphi_1,\cdots,
\varphi_m\in CS_n(A)$ such that $CS_n(A)=\cup_{i=1}^mU(\varphi_i;\frac{\epsilon}{64\lambda_n})$. Since $\hat{\mathbf{S}}$ is bounded, there 
exists an $M_n>0, n\in\mathbb{N}$, such that $\|\hat{a}_n(\varphi)\|\le M_n$ for all $\varphi\in CS_n(A)$ and $a\in S_n$. So $\{\hat{a}_n
(\varphi): \varphi\in CS_n(A), a\in S_n\}$ is totally bounded. The sets $V_i=\{\hat{a}_n(\varphi_i): a\in S_n\}, i=1,2,\cdots,m$, are all 
totally bounded. Thus there are $\frac{\epsilon}{32}$-nets in each of them which are denoted by some sets $\{x_{i1},\cdots,x_{ik_i}\}$ 
where $i=1,2,\cdots,m$. Hence we obtain
\[\cup^m_{i=1}V_i\subseteq\cup^m_{i=1}\cup_{l_i=1}^{k_i}U\left(x_{il_i};\frac{\epsilon}{32}\right)=
\cup^N_{j=1}U\left(z_j;\frac{\epsilon}{32}\right),\]
where $N=k_1+k_2+\cdots+k_m$ and the sets $\{1,\cdots,k_1\}$, $\{k_1+1,\cdots,k_1+k_2\}$, $\cdots$, $\{k_1+\cdots+k_{m-1}+1,\cdots,N\}$ of 
the index $j$ correspond, respectively, to the sets $\{11,\cdots,1k_1\}$, $\{21,\cdots,2k_2\}$, $\cdots$, $\{m1,\cdots,mk_m\}$ of the pair 
of indices $i$ and $l_i$. The points $z_j$ are none other than the points $x_{il_i}$ relabeled in this way.

Let $\Gamma$ denote the finite set of all mappings from $\{1,2,\cdots,m\}$ into $\{1,2,\cdots,N\}$. For any $\gamma\in\Gamma$, we denote
\[O_{\gamma}=\left\{a\in S_n: E^{(n)}_n(\hat{a}_n(\varphi_1),z_{\gamma(1)})<\frac{\epsilon}{32},\cdots,
E^{(n)}_n(\hat{a}_n(\varphi_m),z_{\gamma(m)})<\frac{\epsilon}{32}\right\}.\]
Evidently, we have $S_n=\cup_{\gamma\in\Gamma}O_{\gamma}$. For $\gamma\in\Gamma$ and $a,b\in O_{\gamma}$ and $\varphi\in CS_n(A)$, there exists a 
$\varphi_i$ such that $D_{K_n}(\varphi,\varphi_i)<\frac{\epsilon}{64\lambda_n}$. We obtain
\[\begin{array}{rcl}
&&\|\ll\varphi,a-b\gg\|\\
&=&E^{(n)}_n(\hat{a}_n(\varphi),\hat{b}_n(\varphi))\\
&\le&E^{(n)}_n(\hat{a}_n(\varphi),\hat{a}_n(\varphi_i))+E^{(n)}_n(\hat{a}_n(\varphi_i),z_{\gamma(i)})\\
&&+E^{(n)}_n(z_{\gamma(i)}, \hat{b}_n(\varphi_i))+E^{(n)}_n(\hat{b}_n(\varphi_i),\hat{b}_n(\varphi))\\
&<&\frac{\epsilon}{8}.
\end{array}\]
By the arbitrariness of $\varphi$, we have
\[q_n(\widehat{a-b})=\sup\{\|\ll\varphi,a-b\gg\|: \varphi\in CS_n(A)\}\le\frac{\epsilon}8.\]
In view of Lemma \ref{le:53}, we get $p_n(\widehat{a-b})\le\frac{\epsilon}2$. So $\|a-b\|_n\le\frac{\epsilon}2$ by Lemma \ref{le:52}. 
Hence for any $\gamma\in\Gamma$ there exists an $a_{\gamma}\in S_n$ such that $O_{\gamma}\subseteq U(a_\gamma;\epsilon)$. So
$S_n\subseteq\cup_{\gamma\in\Gamma}U(a_\gamma;\epsilon)$, namely $S_n$ is totally bounded. Since $M_n(A)$ is complete, $S_n$ is norm 
relatively compact, whence (1).
\end{proof}

Now we come to the ``function" level. In \cite{anch}, the point of view on Lip-norms is also nearly the same as the one Kerr has in mind \cite{ker}. In analogy with the 
metric sets of unital $C^\ast$-algebras, we define a {\it metric set} of an operator system $A$, in a strict operator system analogue 
(see Proposition 2.5 and Proposition 4.3 in \cite{ker}), as a subset $K$ of $A$ which 
is norm compact, self-adjoint, and balanced, and convex, and separates the states on $A$. Then we have:

\begin{corollary}\label{co:55}
Let $A$ be an operator system and $K$ a metric set of $A$. For any subset $S$ of $A$ the following conditions are equivalent:
\begin{enumerate}
\item The set $S$ is norm relatively compact.
\item The set of affine functions $\{\hat{a}\in A(\mathcal{S}(A)): a\in S\}$ is bounded and equicontinuous with respect to the 
$w^\ast$-topology on the state space $\mathcal{S}(A)$.
\item The set $S$ is bounded and for every $\epsilon>0$ there exists a $\lambda>0$ such that
\[S\subseteq A_{\epsilon}+\lambda K+\mathbb{C}1,\]
where $A_{\epsilon}=\{a\in A: \|a\|\le\epsilon\}$.
\end{enumerate}
\end{corollary}

\begin{proof}
Since $A$ is complete, $K$ is weakly closed and absolutely convex. So there exists a weakly closed absolutely matrix convex set 
$\mathbf{K}=(K_n)$ with $K_1=K$ and $K_n\subseteq M_n(A)$ for $n\in\mathbb{N}$ (see page 181 in \cite{efwe}). It is easy to verify that $\mathbf{K}$ is a matrix metric set of 
$A$ when view $(A,1)$ as a matrix order unit space.

Fix $k_n\in K_n$ for $n\in\mathbb{N}$. Let $\mathbf{S}=(S_n)$ with $S_1=S$ and $S_n=\{k_n\}$ for $n>1$. By Theorem \ref{th:54}, 
the following conditions are equivalent:

(i) The graded set $\mathbf{S}$ is norm relatively compact.

(ii) The graded set $\hat{\mathbf{S}}=(\hat{S}_n)$, where $\hat{S}_n=\{\hat{a}\in A(\mathcal{CS}(A),M_n): a\in S_n\}$, is bounded and 
equicontinuous with respect to the matrix metric $\mathcal{D}_{\mathcal{K}}$ and the matrix metric $\mathcal{E}^{(n)}$ induced by the matrix norm 
on each $M_n$.

(iii) The graded set $\mathbf{S}$ is bounded and for every $\epsilon>0$ there exists a sequence $\boldsymbol{\lambda}=(\lambda_n)$ with $\lambda_n>0$ such that 
\[\mathbf{S}\subseteq \mathbf{A}_{\epsilon}+\boldsymbol{\lambda}\mathbf{K}+\mathbb{C}\mathbf{1}.\]

Clearly, $\mathbf{S}$ is norm relatively compact if and only if $S$ is norm relatively compact. It is also obvious that (iii) holds 
exactly if $S$ is bounded and for every $\epsilon>0$ there exists a $\lambda>0$ such that 
\[S\subseteq A_{\epsilon}+\lambda K+\mathbb{C}1.\]

Now we need only show that the conditions (2) and (ii) are equivalent. Assume (2). Then $\{\hat{a}\in A(\mathcal{S}(A)): a\in S\}$ is bounded. So there is an $M>0$ such that 
\[q_1(\hat{a})=\sup\{|\hat{a}(\varphi)|: \varphi\in\mathcal{S}(A)\}\le M,\]
for $a\in S$. By Lemma \ref{le:53} and Lemma \ref{le:52}, we have
\[p_1(\hat{a})\le 4q_1(\hat{a})\le 4M,\]
for $a\in S$. Thus $\hat{S}_1=\{\hat{a}\in A(\mathcal{CS}(A),\mathbb{C}): a\in S_1\}$ is bounded. For $n>1$, it is clear that $\hat{S}_n=
\{\hat{k}_n\}$ is bounded by Lemma \ref{le:52}. 

From Proposition \ref{pr:42} and Corollary \ref{co:43}, the $D_{K_1}$-topology on $\mathcal{S}(A)$ 
agrees with the $w^\ast$-topology. Since $\{\hat{a}\in A(\mathcal{S}(A)): a\in S\}$ is equicontinuous with respect to the 
$w^\ast$-topology on $\mathcal{S}(A)$,  $\{\hat{a}\in A(\mathcal{S}(A)): a\in S\}$ is equicontinuous with respect to the $D_{K_1}$-topology. 
Given $\epsilon>0$. For any $\varphi\in\mathcal{S}(A)$, we can find $\delta^\prime=\delta(\epsilon,\varphi,n)>0$ such that 
\[|\hat{a}(\varphi)-\hat{a}(\psi)|<\frac{\epsilon}{4n^2},\]
for all $a\in S$ and $\psi\in\mathcal{S}(A)$ with $D_{k_1}(\varphi,\psi)<\delta^\prime$. That $\mathcal{S}(A)$ is compact in the $D_{K_1}$-topology 
means that there are $\varphi_1,\cdots,\varphi_t\in\mathcal{S}(A)$ such that $\mathcal{S}(A)\subseteq\cup_{i=1}^tU(\varphi_i;
\frac{\delta_i}2)$, where $\delta_i=\delta(\epsilon, \varphi_i,n)$. Take $\delta_0=\frac12\min\{\delta_1,\cdots,\delta_t\}$. Then for 
$\varphi,\psi\in\mathcal{S}(A)$ with $D_{k_1}(\varphi,\psi)<\delta_0$, there is a $\varphi_i$ with $\varphi\in U(\varphi_i;\frac{\delta_i}2)$. 
So $D_{k_1}(\psi,\varphi_i)\le D_{k_1}(\psi,\varphi)+D_{k_1}(\varphi,\varphi_i)<\delta_i$. Hence for any $a\in S$ we have
\[|\hat{a}(\varphi)-\hat{a}(\psi)|\le |\hat{a}(\varphi)-\hat{a}(\varphi_i)|+|\hat{a}(\varphi_i)-\hat{a}(\psi)|<\frac{\epsilon}{2n^2}.\]
So there is a $\delta=\delta(\epsilon,n)>0$ such that 
\[|\hat{a}(\varphi)-\hat{a}(\psi)|<\frac{\epsilon}{2n^2},\]
for $a\in S$ and $\varphi,\psi\in\mathcal{S}(A)$ with $D_{K_1}(\varphi,\psi)<\delta$. 

By Corollary \ref{co:44}, there is a lower semicontinuous matrix Lip-gauge $\mathcal{L}=(L_{n})$ on 
$(A,1)$ such that
\[D_{K_k}(\varphi,\psi)=\sup\{\|\ll\varphi,a\gg-\ll\psi,a\gg\|: L_r(a)\le 1, r\in\mathbb{N}\},\]
for $\varphi,\psi\in CS_k(A)$ and $k\in\mathbb{N}$. From Proposition 3.3 in \cite{wu2}, we have 
\[D_{K_k}(\varphi,\psi)=\sup\{\|\ll\varphi,a\gg-\ll\psi,a\gg\|: a=a^\ast, L_r(a)\le 1, r\in\mathbb{N}\},\]
for $\varphi,\psi\in CS_k(A)$ and $k\in\mathbb{N}$. For $n\in\mathbb{N}$ and $\varphi=[\varphi_{st}], \psi=
[\psi_{st}]\in CS_n(A)$, there are $\phi^{(j)}_{st}\in\mathcal{S}(A)$, $j=1,2,3,4$ and $s,t=1,2,\cdots,n$, such that $\varphi_{st}-\psi_{st}
=\phi^{(1)}_{st}-\phi^{(2)}_{st}+i(\phi^{(3)}_{st}-\phi^{(4)}_{st})$ by Lemma 4.1 in \cite{wu2}. For $b=b^\ast\in M_r(A)$ with 
$L_r(b)\le 1$, we have
\[\begin{array}{rcl}
\left\|\ll\phi^{(1)}_{st},b\gg-\ll\phi^{(2)}_{st},b\gg\right\|&\le&\|\ll\varphi_{st},b\gg-\ll\psi_{st},b\gg\|\\
&=&\|(e_s\otimes1_r)(\ll\varphi,b\gg-\ll\psi,b\gg)(e^\ast_t\otimes 1_r)\|\\
&\le&\|\ll\varphi,b\gg-\ll\psi,b\gg\|\\
&\le&D_{K_n}(\varphi,\psi),\end{array}\]
where $e_s=[0 \ \cdots \ 0 \ 1_s \ 0 \ \cdots \ 0]$. So we have $D_{K_1}(\phi^{(1)}_{st}, \phi^{(2)}_{st})
\le D_{K_n}(\varphi,\psi)$ for $s,t=1,2,\cdots,n$. Similarly we have $D_{K_1}(\phi^{(3)}_{st}, \phi^{(4)}_{st})
\le D_{K_n}(\varphi,\psi)$ for $s,t=1,2,\cdots,n$. When $D_{K_n}(\varphi,\psi)<\delta$, we have $D_{K_1}
(\phi^{(1)}_{st},\phi^{(2)}_{st})<\delta$ and $D_{K_1}(\phi^{(3)}_{st},\phi^{(4)}_{st})<\delta$, and hence 
\[\begin{array}{rcl}
E^{(1)}_n(\hat{a}_n(\varphi),\hat{a}_n(\psi))&=&\|\varphi(a)-\psi(a)\|\\
&\le&\sum_{s,t=1}^n|\varphi_{st}(a)-\psi_{st}(a)|\\
&\le&\sum_{s,t=1}^n\left(\left|\phi^{(1)}_{st}(a)-\phi^{(2)}_{st}(a)\right|+\left|\phi^{(3)}_{st}(a)+\phi^{(4)}_{st}(a)\right|\right)<\epsilon,\end{array}\]
for $a\in S$. Hence $\hat{S}_1$ is equicontinuous with respect to $\mathcal{D}_{\mathcal{K}}$ and $\mathcal{E}^{(1)}$. For $n>1$, $\hat{S}_n
=\{\hat{k}_n\}$ is equicontinuous with respect to $\mathcal{D}_{\mathcal{K}}$ and $\mathcal{E}^{(n)}$ since $\hat{k}_n\in\mathcal{C}(\mathcal{CS}(A), M_n)$, 
whence (ii). The implication (ii) implies (2) is obvious. So conditions (2) and (ii) are equivalent.
\end{proof}

\begin{remark}
(1) According to Definition 3.1 in \cite{anch}, a metric set in a unital $C^\ast$-algebra may be not self-adjoint. For example, let 
$G=\{e, g_1, g_1^{-1}, g_2, g_3, \cdots\}$ be a countable group, where $e$ is the identity element of $G$ and $g_1^2\neq e$. We denote 
by $K$ the closed convex hull of the set $S=\{\delta\lambda_e: \delta\in\mathbb{C}, |\delta|\le 1\}\cup\{\theta\lambda_{g_1}: \theta\in\mathbb{C}, 
|\theta|\le\frac12\}\cup\{\alpha\lambda_{g_n}+\beta\lambda^\ast_{g_n}: |\alpha|+|\beta|\le\frac1{n+1}, \alpha,\beta\in\mathbb{C}, 
n\ge 2, n\in\mathbb{N}\}$. Then $K$ is norm compact, balanced and convex, and separates the states on $C^\ast_r(G)$. But $K$ is not self-adjoint.

(2) We would like to point out that the definition of a metric set in \cite{anch} should contain the condition of self-adjointness. One 
reason is that a Lipschitz seminorm $L$ for a $C^\ast$-algebra $\mathcal{A}$ should satisfy $L(a^\ast)=L(a)$ for $a\in\mathcal{A}$ (see 
page 6 in \cite{ri4} or Proposition 2.5 and Proposition 4.3 in \cite{ker}). Another reason is that without the self-adjointness, we can not 
get
\[\forall h\in\mathcal{H}\ \forall\gamma\in\mathcal{A}_2^\ast\cap\{\mathbb{C}I\}^{\bot}\cap\delta(\mathcal{K}^{\circ}): 
|\gamma(h)|\le 2\epsilon\]
from
\[\forall h\in\mathcal{H}\ \forall\gamma\in(\mathcal{A}_h^\ast)_2\cap\{\mathbb{C}I\}^{\bot}\cap\delta(\mathcal{K}^{\circ}): 
|\gamma(h)|\le \epsilon,\]
(see page 258 in \cite{anch}) because it is not guaranteed that $\frac12(\gamma+\gamma^\ast)$ and $\frac12(\gamma-\gamma^\ast)$ belong to 
$\delta(\mathcal{K}^\circ)$ for $\gamma\in\delta(\mathcal{K}^\circ)$.
\end{remark}

\section*{Acknowledgements}

I would like to thank Marc Rieffel for valuable discussions and suggestions. This research was partially supported by 
Shanghai Priority Academic Discipline, China Scholarship Council and National Natural Science Foundation of China.

\bibliographystyle{amsplain}

\end{document}